\input{BoxedEPS.tex}
\SetepsfEPSFSpecial
\HideDisplacementBoxes
\documentclass[11pt]{amsart}
\begin{document}
\title[]{The Modular Tree of Pythagorus}
\author[]{Roger C. Alperin}

\newtheoremstyle{citing}
  {5pt}
  {5pt}
  {\itshape}
  {}
  {\bfseries}
  {.}
  {.5em}
  {}

\theoremstyle{citing}
\addtolength{\textwidth}{100pt}
\hfuzz 100pt
\newtheorem{theorem}{Theorem}[section]
\newtheorem{proposition}[theorem]{Proposition}
\newtheorem{lemma}[theorem]{Lemma}
\newtheorem{corollary}[theorem]{Corollary}
\newtheorem{remark}[section]{Remark}
\newtheorem{definition}[section]{Definition}
\newtheorem{question}{Question}

\def\cal{\mathcal}

\def\proof{{\bf {\noindent}Proof. }}
\newsymbol\bsq 1004
\def\endproof{\hfill${\bsq}$\bigskip}

\def\question{\par{\bf {\smallskip}{\noindent}Question: }}
\def\definition{\par{\bf {\smallskip}{\noindent}Definition: }}
\def\remark{\par{\bf {\smallskip}{\noindent}Remark: }}
\def\note{\par{\bf {\smallskip}{\noindent}Note: }}
\def\example{\par{\bf {\smallskip}{\noindent}Example: }}

\def\bar{\overline}

\newcommand\R{\mbox{\bf R}}
\newcommand\C{\mbox{\bf C
}}
\newcommand\Z{\mbox{\bf Z}}
\newcommand\Q{\mbox{\bf Q}}

\newcommand\SL{{\rm SL}}
\newcommand\PSL{{\rm PSL}}
\newcommand\GL{{\rm GL}}
\newcommand\M{{\rm M}}
\newcommand\E{{\rm E}}

\def\psmod{\rm PSL_2(\bf Z)}
\def\psmodR{\rm PSL_2(\bf R)}
\def\smod{\rm SL_2(\bf Z)}
\def\smodR{\rm SL_2(\bf R)}
\def\pmod{\rm PGL_2(\bf Z)}
\def\mod{\rm GL_2(\bf Z)}
\def\two{{\bf Z}^2}
\def\mat{\rm M_2(\bf Z)}
\def\matR{\rm M_2(\bf R)}
\maketitle
\centerline{Department of Mathematics and Computer Science}
\centerline{ San Jose State University, San Jose, CA 95192 USA}

\section{Introduction}
The Pythagorean triples of integers satisfying $x^2+y^2=z^2$ have been studied and enumerated since Babylonian times. Since Diophantus, it has been
known that this set of triples is related to the standard rational parameterization of the unit circle, $(\frac{t^2-1}{t^2+1},\frac{2t}{t^2+1})$. The
Pythagorean triple solutions, which are relatively prime, have the elementary and beautiful characterization as integers
$x=m^2-n^2, y=2mn, z=m^2+n^2$ for relatively prime integers $m, n$.  One can also realize that the Pythagorean triples are related to the Gaussian
integers, $\Z[i]$, the lattice in the complex numbers with integer coordinates, $u=x+iy$, where $x, y$, are integers. The Pythagorean equation 
$N(u)=u\bar{u}=(x+iy)(x-iy)=z^2$ is an equation among Gaussians.  It is perhaps, not surprising, that the Pythagorean triples are just the
squares of the set of Gaussian integers, that is
$v=(m+ni)^2=(m^2-n^2)+(2mn)i$ gives the endpoint of the hypotenuse of a right triangle with integer length sides. This Gaussian square $v$ obviously has
a norm, $N(v)$, which is a square.
Both of these expressions secretly involve the double angle formulae for sine and cosine since the stereographic projection formula uses the
central and chordal angles.

By considering the action of the modular group, $\Gamma=\psmod$, by conjugation on the
set of all integer matrices
$\mat$, we can blend together these two perspectives.  We shall show (Corollary \ref{Gamma})
that the Pythagorean triples can be identified with an orbit of a subgroup (of index 6),
$\Gamma(2)$, which is generated by  the images of $U^2=\begin{pmatrix} 1&2\\0&1\end{pmatrix}$ and
$L^2=\begin{pmatrix} 1&0\\2&1\end{pmatrix}$. Moreover 
$\Gamma(2)$ is the free product of the subgroups generated by $U^2$ and $L^2$ (Proposition
\ref{freeprod}). The consequence of the free group structure for $\Gamma(2)$ is that the Pythagorean
triples can be enumerated  by words, for each $k\ge 0$, $U^{\pm 2 n_k}L^{\pm 2
m_k}U^{\pm 2 n_{k-1}}\cdots L^{\pm 2 m_1}U^{\pm 2 n_1}L^{2 m_0}$, for integers $m_i\ne 0$, $n_i\ne
0$, $0\le i\le k$.

\begin{figure}[t]
\label{PT}
\centerline{\BoxedEPSF{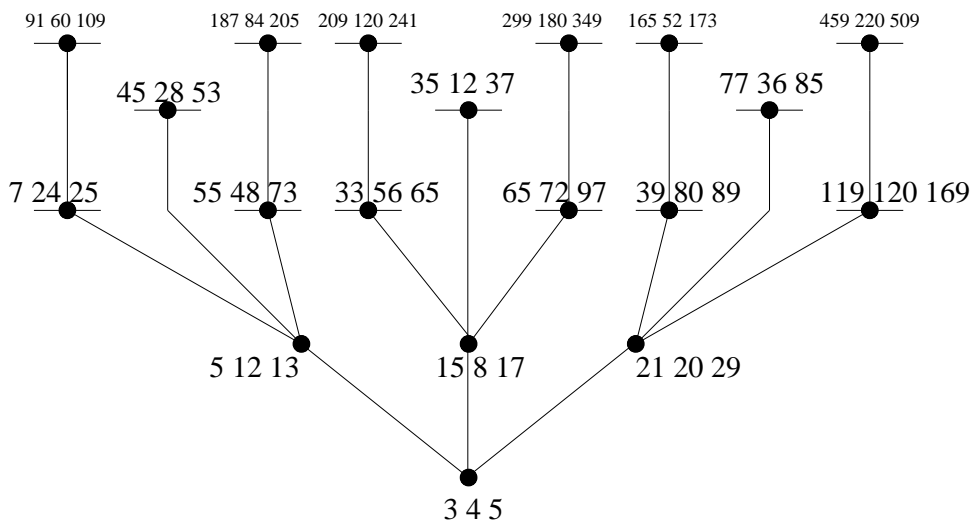 scaled 700}}
\caption{\it Tree of Pythagorus}
\end{figure}

Using a breadth-first enumeration of these group elements we obtain the following curious
result (see Figure \ref{PT}). 

\begin{theorem}
\label{tree} The set of primitive Pythagorean triples has the structure of a complete infinite rooted ternary tree.
\end{theorem}
\bigskip

\section{Conjugation Action and Similarity}
\begin{proposition}
\label{cone}
Any integer matrix $X$, satisfying  $X^2=0$ is
$X=\begin{pmatrix} x&y\\z&-x\end{pmatrix}$ for integers $x, y, z$ satisfying $yz+x^2=0$ and conversely.
\end{proposition} 
\proof  Suppose the integer matrix $X$, satisfies
$X^2=0$. Since the matrix $X$ also satisfies its characteristic polynomial, it follows that $tr(X)X=det(X)I$. Since it is impossible for
a non-zero $X$ to be a non-zero multiple of the identity matrix, it follows that $tr(X)=det(X)=0$, and the conclusion of the lemma is immediate. 
\endproof

The group $\mod$ of invertible linear transformations of
$\two$ acts via conjugation on all two by two matrices,
$\mat$, by 
$X\rightarrow TXT^{-1}$.  This action preserves the trace and determinant and hence also the  `cone' of nilpotent matrices $\cal{N}_2=\{ X\in\ $ $\mat$ $\ |\
X^2=0\}$ as described by Proposition \ref{cone}. 

The similarity class of the matrix 
$X$ is $[X]=\{ TXT^{-1} |\ T\in \mod\}$. 
Let $E=\begin{pmatrix} 0&-1\\0&\ 0\end{pmatrix}$;  Notice  that $E^t\in [E]$ using conjugation by  $\begin{pmatrix}
0&1\\1&0\end{pmatrix}$; also $-E\in [E]$ using conjugation by the matrix  $D=\begin{pmatrix} -1&0\\0&1\end{pmatrix}$. 

Let $\bf N$ denote the set of  non-negative integers. Let ${\cal E}_\lambda= [\lambda E]$, $\lambda\in \bf N$; it is clear that ${\cal
E}_\lambda=\lambda[E]$.  Also, it is an easy check to see that ${\cal E}_\lambda$ and ${\cal E}_\mu$, for $\lambda\ne \mu$, are disjoint.

\begin{proposition}
A matrix $X\in {\cal N}_2$ is similar to $\lambda E$, for a unique  $\lambda\in \bf N$. Thus  ${\cal N}_2=\bigcup_{\lambda\in \bf {N}}
\cal E_\lambda$ is  the  disjoint union of similarity classes.
\end{proposition} 
\proof If any entry of the matrix $X$ is zero, then it follows immediately from Proposition \ref{cone} that $X$ is either the
 zero matrix, a multiple of $E$, or a multiple of  $E^t$. 
 Thus, in this case $X$ is similar to either the zero matrix or $E$.

Suppose now that none of the entries of the matrix $X$ is zero; say $X=\begin{pmatrix}
x&y\\z&-x\end{pmatrix}$.  Since the determinant is zero, the non-zero rows are dependent over the rationals; furthermore, after factoring
out the greatest common divisor, say $\lambda$, of the entries of $X$, we find that there are 
relatively prime integers $m$ and $n$ so that $mx=nz$ and $my=-nx$.

Now, using the relatively prime condition as above,  we obtain the divisibility:  $m|z$, $m|x$, $n|x$, $n|y$. Rewriting then $x=mnx_1$, $y=ny_1$, $z=mz_1$;
canceling where possible then from the original equations, it  follows that $mx_1=z_1$, $-nx_1=y_1$; thus $m|z_1$ and $n|y_1$. Hence, 
$x=mn\lambda, y=-n^2\lambda, z=m^2\lambda$, $\lambda\in \Z$.

Thus the original matrix is a multiple of 
 a matrix, which has a factorization,   $$\begin{pmatrix} mn&-n^2\\m^2&-mn\end{pmatrix}=\begin{pmatrix}
n\\m\end{pmatrix}\begin{pmatrix}
m&-n\end{pmatrix}.$$  
Consequently, the action of conjugation is just a mixture of the usual action on $\two$ as row and
as column vectors:
$$TXT^{-1}=T\begin{pmatrix}
n\\m\end{pmatrix}\begin{pmatrix}
m&-n\end{pmatrix}T^{-1}.$$
The conjugation action of 
$T=\begin{pmatrix} a&b\\c&d\end{pmatrix}$ gives $$TXT^{-1}=\pm\begin{pmatrix}
(an+bm)(cn+dm)&-(an+bm)^2\\(cn+dm)^2&-(an+bm)(cn+dm)\end{pmatrix}.$$ 

Now it is easy to see, using the Euclidean algorithm, that for any relatively prime integers $m, n$ there is an $T\in \smod$, the elements of $\mod$ of
determinant one, so that
$T\begin{pmatrix} n\\m\end{pmatrix}=\begin{pmatrix}
1\\0\end{pmatrix}$. Hence using the conjugation of such a $T$ on $X$ as described above gives the matrix $E$.
Thus any non-zero integer matrix $X$ with relatively prime entries is similar to $E$. 
\endproof

\section{Enumerating Pythagorean Triples}

Consider the matrices in $\mat$ which are similar to $E$; there are relatively
prime integers $m,n$ so that the matrix is  $$\begin{pmatrix} mn&-n^2\\m^2&-mn\end{pmatrix}= 
\frac{1}{2}\begin{pmatrix} C&S-N\\S+N&-C\end{pmatrix}.$$ The integers $S, C,
N$  give a Pythagorean triple satisfying $S^2+C^2=N^2$.

We now have  a
method of generating Pythagorean triples. Consider the conjugation action of $\mod$ on $E$ and
determine the orbit. Of course this is just the same as the set of cosets $\mod$ by the stabilizer
of $E$.   The stabilizer in $\mod$ of the matrix
$E$ is easily seen to be the subgroup generated by 
$U$ and the scalar matrix $-I$. However, the situation is a bit subtle. These cosets parameterize the different matrices on the left side of the equation above,
but not the different Pythagorean triples, which are obtained from the right hand side. For primitive ($m$, $n$ relatively prime) Pythagorean triples
$S$, $C$,
$N$ we only want to count the really different solutions; fix $m$ to be even and $n$ to be odd and assume both $m$ and $n$ are positive; this specifies 
$C$ as even and positive, $N$ is positive, and $S$ is odd. 
To enumerate the distinct Pythagorean triples as described above, we shall use matrices of determinant 1 with even $2-1$ entry, $m$; it
follows that the $1-1$ entry $n$ is odd since the matrix $T$ is invertible; hence also $m$ and
$n$ are relatively prime; thus, 
$TET^{-1}=\begin{pmatrix}
mn&-n^2\\m^2&-mn\end{pmatrix}$. Since the  effect of the matrix $-I$, the
negative identity, (on triples $S$, $C$, $N$) is trivial,  we  mod out that action. So we
have either both $m$ and $n$ are positive or opposite signs; since both situations occur the
counting the cosets will double count the Pythagorean triples. We want to use the cosets to
enumerate the Pythagorean triples in a 1-1 fashion; and so we look for a way avoid this
duplication. We shall develop the group theory a bit more to aid in this.

The homomorphism $\tau: \rm{PSL_2({\bf Z})}\longrightarrow \rm{PSL_2({\bf Z}_2)}$, induced by
reducing mod 2,  has  kernel denoted, $\Gamma(2)$.  Since the  image of $\tau$ is the group
$\rm{PSL_2({\bf Z}_2)}$ (it is a non-abelian group of order 6) the group 
$\Gamma(2)$ is of index 6 in $\psmod$ and  normal. It is easy to see that $\Gamma(2)$ is generated
by (the images of) $U^2=\begin{pmatrix} 1&2\\0&1\end{pmatrix}$ and $L^2=\begin{pmatrix}
1&0\\2&1\end{pmatrix}$. Moreover, these elements generate the group as a free product. This will
be  useful for the enumeration of Pythagorean triples. For any relatively prime integers $n$, $m$,
one odd and one even,  there is a matrix in $\psmod$ with those as the entries of the first
column; but multiplying by an appropriate coset representative, we have that these elements $n$,
$m$ (up to sign) are the first column of a matrix in $\Gamma(2)$. 

\begin{proposition}
\label{freeprod} $\Gamma(2)$ is the free product of the infinite cyclic subgroups generated by  $U^2$ and $V^2$.
\end{proposition}
\proof This result is easily proven using the  structure of $PSL_2({\bf Z})$ as the free product
\cite{A1}, of the subgroups of orders 2 and 3 generated  by 
$A=\begin{pmatrix} 0&1\\-1&0\end{pmatrix}$ and $B=\pmatrix 0&-1\\1&1 \endpmatrix$. Since
$L=AB$, $L^2=ABAB$, $U=AB^{-1}$, $U^2=AB^{-1}AB^{-1}$, then  any alternating word in $U^{\pm 2}$ or $V^{\pm 2}$ is also 
alternating in
$A, B^{\pm 1}$ and hence non-trivial. 
\endproof

The enumeration of triples is based on the enumeration of cosets of $\Gamma(2)$ modulo the stabilizer of $E$, which is the subgroup generated by $U^2$.  Because of
the free product structure the coset representatives are just the words  in
$L^2$, $U^2$; since we are considering cosets of the subgroup generated by $U^2$, the rightmost letter of a coset
representative is $L^{\pm 2}$.  We do a breadth first enumeration, adding $U^{\pm 2}$ or $L^{\pm 2}$ on the left of a string (alternating); for example, if the last
group element on the left  is
$L^2$ we can add $L^2$,
$U^2$ or
$U^{-2}$.

For any matrix $M\in \mat$, let $\delta(M)$ be the matrix with the same off diagonal entries as
$M$ and with signs changed on the diagonal. We want to
avoid counting both $M$ and $\delta(M)$ in the enumeration of cosets of $\Gamma(2)$ since they produce the same Pythagorean triple from the action on 
$E$. The conjugation action of the matrix
$D$ and  multiplication by $-I$ has the same effect as $\delta$,  $\delta(M)=-(DMD^{-1})=DMD^{-1}$ (projectively).
Thus $\delta$ gives an outer automorphism of order 2 of $\Gamma(2)$;  conjugation by $D$.
On the generators of $\Gamma(2)$ the automorphism is: $$\delta(U^2)=U^{-2},\delta(L^2)=L^{-2}.$$ 

We then can effectively describe the distinct
coset representatives which will enumerate the distinct Pythagorean triples  using the alternating property for the words of the free group $\Gamma(2)$ in the
generators $L^2$ and $U^2$.  Also, using
$\delta$ we see that to obtain distinct triples we can initialize with 
 ${\cal L}^+_0=\{ L^2\}$, ${\cal L}^-_0=\{\}$, ${\cal U}_0^\pm=\{\}$ and define inductively, $i\ge 0$,
\begin{align*}
{\cal L}^+_{i+1}&=\{\ L^2X\ |\ X\in {\cal L}^+_i\cup {\cal U}^\pm_i\ \},\\ 
{\cal L}^-_{i+1}&=\{\ L^{-2}X\ |\ X\in {\cal L}^-_i\cup {\cal U}^\pm_i\ \},\\ 
{\cal U}^+_{i+1}&=\{\ U^2X\ |\ X\in {\cal U}^+_i\cup {\cal L}^\pm_i\ \},\\
{\cal U}^-_{i+1}&=\{\ U^{-2}X\ |\ X\in {\cal U}^-_i\cup {\cal L}^\pm_i\ \}.
\end{align*}

Thus we have an effective enumeration of non-trivial primitive Pythagorean triples using the
coset representatives in the disjoint union
$$\cal P=\bigcup_{i\ge 0} ({\cal L}^{\pm}_i\cup {\cal M}^\pm_i).$$ 
The recurrence relation which counts the number of elements $p_k$ in $\cal P$  of level $k$  is  $p_{k}=3p_{k-1}$, since we can add any one of 3 elements to the
left of a given alternating string to keep it alternating. We use $\Pi$ to denote the triples obtained after conjugation of $E$ by the elements of $\cal P$. 

For each $T\in \Pi$, the entries determine a unique primitive Pythagorean triple of positive integers as
follows:
$$|S|=|T_{2,1}+T_{1,2}|,|C|=|2T_{1,1}|,|N|=|T_{2,1}-T_{1,2}|.$$

\begin{theorem}
\label{Gamma}
The primitive Pythagorean triples  are in 1-1 correspondence with the set
$\Pi$ of certain of the
the coset representatives  of the subgroup generated by $U^2$ in $\Gamma(2)$. This set $\Pi$ can be enumerated as a union of triples of level $k$
indexed by subsets $\Pi_k$ of size $3^k$.
\end{theorem}

To make the tree we  connect the $j^{th}$ element of level $i$ with the 3 elements of level $i+1$ numbered $3(j-1)+1,3(j-1)+2,3(j-1)+3$ corresponding to the 3
possible elements used to make the next level by the conjugation action of three of $L^{-2}, L^2,  U^{-2}, U^2$; this  yields triples, up
to sign, respectively, which can further be exprssed in terms of the coordinates
$S=m^2-n^2,C=2mn,N=m^2+n^2$,
\begin{align*}
{\bf L_-}&=[m^2 -4 m n   + 3 n^2 , 2m n - 4 n^2,   m^2-4 m n +   5 n^2]\\
{\bf L_+}&=[m^2  + 4 m n +  3 n^2 , 2m n + 4 n^2,   m^2+4 m n +  5 n^2]\\
{\bf U_-}&=  [-n^2  +4 m n -  3 m^2, 2m n - 4 m^2 , n^2 -4 m n  +5 m^2 ]\\
{\bf U_+}&=  [ -n^2  -4 m n - 3 m^2, 2m n + 4 m^2  , n^2 + 4 m n + 5 m^2 ].
\end{align*}

Amazingly, we only have to compute the three transformations $\bf L_{+}, \bf U_{-}, \bf U_{+}$ when the coordinates
$S, C, N$ are all positive to get the next level in the tree. When $S, C, N$ are all
positive then the previous triangle is 
$\bf L_{-}$. This is easy to see, since in this case $m> n\ge 0$,  and therefore $2mn>2mn-4n^2
>-2mn$ so that the second coordinate of 
$\bf L_{-}$, $2mn-4n^2$, is smaller in absolute value than $C=2mn$. Similarly, the second
coordinate of each of $\bf L_{+}, \bf U_{-}, \bf U_{+}$ is greater than
$C$. Furthermore, with variations on this kind of argument, we can modify these three
transformations to obtain three new transformations, so that all the coordinates are positive,
and so preserve the positive quadrant of the cone $x^2+y^2=z^2$.

\begin{align*}
{\cal L_+}&=[m^2  + 4 m n +  3 n^2 , 2m n + 4 n^2,   m^2+4 m n +  5 n^2]\\&=[S,-C,N]+ 2(N-S+C) [ 1 , 1, 1]\\
{\cal U_-}&=  [n^2  -4 m n +  3 m^2,  4 m^2- 2mn , n^2 -4 m n  +5 m^2 ]\\&=  [-S,C,N]+ 2(N+S-C) [1 , 1 , 1 ]\\
{\cal U_+}&=  [ n^2 + 4 m n + 3 m^2, 2m n + 4 m^2  , n^2 + 4 m n + 5 m^2 ]\\&= [-S,-C,N]+ 2(N+S+C) [1  , 1  , 1 ].
\end{align*}

One can easily check that each of the differences between coordinates in
$[m^2-n^2,2mn,m^2+n^2]$, namely, $(m-n)^2, 2n^2, 2mn-m^2+n^2$, persists to the next
level (except for a sign change on the last one) for exactly one of ${\cal U_-}, {\cal L_+}, {\cal
U_+}$. These three new points are so represented in the tree from left to right (see Figure
\ref{PT}) but connected to the level
$i$ point. For example one can follow the sequence with difference 1 between  the  $z$ and $y$
coordinates, in the tree as the path $$[3,4,5],[5,12,13],[7,24,25],[9,40,41], \cdots .$$ The
difference
$1$  between the $y$ and $x$ coordinates however will alternate, $$[3,4,5], [21,20,29],
[119,120,169], [697,696,985], \cdots.$$ In fact any difference at level $i$ will continue up the
levels of the tree in a similar fashion. The difference $-7$ occurs (alternating sign) for
$$[15,8,17], [65,72,97], [403,396,565], \cdots.$$

One can ask for a  characterization of all those integers which can be differences. 
This is just a question about the representations, for $m$ even and $n$ odd, of integers by
the forms
$P=(m-n)^2$,
$Q=2n^2$,
$R=(n+m)^2-2m^2$ or $-R=(m-n)^2-2n^2$. The first can take on any odd square value, the second, any
twice an odd square, while the third case is more restrictive. The norm form $X^2-2Y^2=D$ can be
solved for 
$Y$ and  $X$ for exactly those integers $D$ which are norms in the ring of integers for the field
$\Q(\sqrt{2})$. However we require that $Y$ is even, so this
is more  restrictive; since $Y$ is even, consider $x^2-8y^2=D$; then it also follows that any
solution when $D$ is odd has $x$ even. 
Consider this equation modulo $p$ for any prime divisor of $D$ then it follows that  2 is a
quadratic residue modulo $p$.  Using quadratic reciprocity, it follows  that every  prime divisor
of $D$  must be $\pm 1\ mod\ 8$. However, for a prime divisor of $D$, $p=-1\ mod\ 8$, it follows that
$x^2=p=-1\ mod\ 8$; which however is impossible. It follows that any odd prime factor of $D$ having
odd exponent must be congruent to 1 modulo 8. Since the ring of integers in $\Q(\sqrt{2})$ is a
unique factorization domain it follows that for every $p$ dividing $D$  there is an
element of norm $p$ when $p=1\ mod\ 8$ (not just a multiple). Multiplying these together gives a
solution for $D$. Furthermore, for the minimal  positive representation of $D$ by this form $R$ we
obtain an infinite path in the tree, since we can solve for $m$ and $n$ from the equations $X=n+m,
Y=m$.
 
Given the absolute values of $x, y$, solutions to $x^2-2y^2=R>0$, then we can solve
for  values of $m, n$ yielding a Pythagorean triple. Choose $n+m=x>0$ and $m=y<0$ then
$n=x-m>0$. Hence the Pythagorean triple of positive integers $[n^2-m^2,-2mn,m^2+n^2]$ has first
difference  $R=n^2-m^2-(-2mn)$. For example to find the first difference 17, one has the root
Pythagorean triple  coming from the solution to
$5^2-2\cdot 2^2=17$; giving rise to the solution, $5+2\sqrt{2}$, so $m=-2,
n+m=5$,  yields the triple
$S=45, C=28, N=53$; to solve for $D=-17$, we can multiply by the smallest unit of negative norm,
$-1+\sqrt{2}$, obtain the solution $(-1+\sqrt{2})(5+2\sqrt{2})=-1+3\sqrt{2}$ so
solving $-R$, we obtain $n=3, m-n=1$, which yields $[7,24,25]$.
One can also use a method of Lagrange as described in \cite{C}, to solve $x^2-8y^2=p$. For example
when $p=89$, one finds the smallest $z_0<89/2$, so that
$(z_0^2-8)/89$, is an integer, say $q_1$; here
$z_0=39$, $q_1=17$; then solve for the smallest $z_1<17/2$, with $(z_1^2-8)/17=q_2$. In this case
$z_1=5$, $q_2=1<\sqrt{8}$ and the algorithm terminates. It now follows that $(39- 2\sqrt{2})/(5+
2\sqrt{2})$ has norm 89. This gives rise to the solution $11+ 4\sqrt{2}$, an element of norm 89;
the number $3-7\sqrt{2}$ has norm -89.
The solution  $11+4\sqrt{2}$, yields $m=-4,
n+m=11$, so $n=15$  yields the triple
$S=209, C=120, N=241$; from the solution $3-7\sqrt{2}$  
to $-R$, we obtain $n=7, m-n=3$,  which yields the triple $[51,140, 149]$. The last difference is
9 which occurs in   Figure \ref{PT} at $[33,56,65]$. As a check one may compute the three
neighbors of this point yielding $[209, 120, 241]$, $[51, 140, 149]$, $[275, 252, 373]$.
Thus $[51,140, 149]$ is the root of those points with first difference -89.

\par
\email{\tt  alperin@mathcs.sjsu.edu}
\medskip\par
AMS Classification: 11D09
\end{document}